\documentclass[12pt,a4,reqno]{amsart}
\usepackage{amsmath,amsfonts,amssymb}
\usepackage{version,xspace}
\usepackage{color,marginnote}
\usepackage{enumerate}
\usepackage{kurier}
\usepackage{hyperref}

\usepackage{marginnote}
\usepackage[normalem]{ulem}

\newtheorem{Theorem}{Theorem}[section]
\newtheorem{Ex}[Theorem]{Example}
\newtheorem{Def}[Theorem]{Definition}
\newtheorem{Lemma}[Theorem]{Lemma}
\newtheorem{Proposition}[Theorem]{Proposition}
\newtheorem{Corollary}[Theorem]{Corollary}

\newtheorem{kNote}[Theorem]{Note}
\newtheorem{Exer}[Theorem]{Exercise}
\newtheorem{Rem}[Theorem]{Remark}
\newtheorem{Remm}[Theorem]{Remarks}
\newtheorem{kNotes}[Theorem]{Notes}
\newenvironment{Example}{\begin{Ex}\normalfont}{\end{Ex}}
\renewcommand\rm\normalfont
\renewcommand\bf\bfseries
\renewcommand\it\itshape
\newenvironment{Definition}{\begin{Def}\normalfont}{\end{Def}}

\newenvironment{Proof}{\proofn}{\hfill$\square$\par\medskip}

\newenvironment{Remarks}{\begin{Remm}\rm}{\end{Remm}}

\newcommand{\tp}[3]{\{#1,#2,#3\}}
\def\D#1#2{#1 \mbox{\footnotesize $\,\square\,$} #2}
\newcommand\q[1]{Q_{#1}}
\def\GL#1{\mbox{\sffamily GL}(#1)}
\newcommand\g[2]{g_{#1}(#2)}
\newcommand\B[1]{B_{#1}}

\newcommand\JBT{$J\kern-1.5pt B^*$-triple}
\newcommand\C{\mathbb{C}}
\newcommand\R{\mathbb{R}}
\newcommand\N{\mathbb{N}}
\newcommand\T{\mathbb{T}}

\newcommand\df[1]{{\emph{#1}}}
\newcommand\Id{{I}}
\renewcommand\L{\mathcal{L}}
\renewcommand\H{\mathcal{H}}
\newcommand\cch[1]{\overline{\operatorname{co}}#1}
\newcommand\Isom{\operatorname{Isom}}
\newcommand\rank{\operatorname{rank}}

\def\ol{\overline}

\def\ni{\noindent}

\newcommand\half{{\frac 12}}
\newcommand\Aut[1]{{\mathsf{Aut}}\kern1pt(#1)}
\def\inv{^{-1}}

\newcommand\absval[1]{\left\vert #1 \right\vert}

\newcommand{\ip}[2]{\langle #1, #2 \rangle}

\renewcommand\tilde{\widetilde}

\newcommand{\norm}[1]{\lVert #1 \rVert}
\newcommand\del{\partial}
\newcommand\proofn{\ni{\bf Proof. \hskip2mm}}

 \newcommand\underlabel[2]{\underset{#1}{\underbrace{#2}}}

\reversemarginpar

\author{M. Mackey \and\ P. Mellon\\  University College Dublin}
\email{macǩey@maths.ucd.ie, pmellon@maths.ucd.ie}

\date{\today} 
\title[Bergmann-Shilov Boundary]{The Bergmann-Shilov boundary of a bounded
  symmetric domain}
\begin{document}
\maketitle

\begin{abstract} We show that there are many sets in the boundary of a
  bounded symmetric domain that determine the values and norm of
  holomorphic functions on the domain having continuous extensions to the boundary. We provide an analogue of the
  Bergmann-Shilov boundary for finite rank \JBT s.
 \end{abstract}

\section*{Introduction}

Recall that an open unit ball $B$ of a complex Banach space $Z$ is
homogeneous with respect to biholomorphic mappings if, and only if,
$Z$ carries an algebraic structure that renders it a \df{\JBT},
defined below \cite{Kaup_RMT}.  \JBT{}s include Hilbert spaces,
$C^*$-algebras and the classical Hermitian symmetric spaces known as
Cartan factors. Bounded symmetric domains are the infinite
dimensional analogues of the Hermitian symmetric spaces but, by Kaup's
Riemann Mapping Theorem \cite{Kaup_RMT}, we may alternatively
introduce them as those domains
in a Banach space which are biholomorphically equivalent to the unit
ball of a \JBT.   This then is a natural category of Banach spaces in
which to study holomorphic functions and demonstrates a rich interplay of complex
and functional analysis, geometry and non-associative algebraic
structures.

The first section of this paper consititutes a very brief introduction to
the basic facts about \JBT s, including the concrete description of
the biholomorphic mappings (automorphisms) of the unit ball in terms
of the triple product.   Tripotents and their order structure are
introduced and linked to the extreme points of the unit ball.   In
Section 2 we prove that the automorphisms (which naturally extend to the
boundary) preserve the set of extreme points as well as the unitary
tripotents.   

These facts allow us in Section 3 to express the closed
unit ball as the closed convex hull of the set of maximal and unitary
tripotents, when these sets are non-empty and, more generally, as the closed convex hull of the orbit 
of any $v \in \del B$ under $\Aut B$, denoted $G_v$, or its orbit under the connected identity component of $\Aut B^0$, denoted $G_v^0$ . 

In Section 4 we see that any holomorphic function on $B$ which
extends continuously to the boundary has its values inside the ball, and hence
its norm, determined by its values on certain subsets of the boundary, which we call determining sets.
We show that such determining sets include the set of extreme points and the set of  unitary tripotents (if these exist) and, more generally, the orbits $G_v$ and $G_v^0$, of any $v \in \del B$. We  recall
the notion and structure of boundary components of a finite rank
\JBT. We  prove that automorphisms of $B$ map holomorphic boundary components onto holomorphic boundary components and, more importantly, they preserve the rank of such boundary components.

Section 5 culminates in our main result, namely,
that the set of extreme points of a finite rank triple
acts analagously to the Bergmann-Shilov boundary in finite dimensions
by providing the smallest closed subset of $\ol B$ which
determines the norm of all holomorphic functions of the ball which have a
continuous extension to the boundary.

\section{Notation and background}

Throughout, $\Delta=\{z\in\C:|z|<1\}$ and $\T=\del\Delta$. For $X$ and $Y$ complex Banach
spaces, $\L(X,Y)$ denotes the space of continuous linear maps from $X$
to $Y$,\ $X^\prime=\L(X,\C)$ and $\L(X)=\L(X,X)$. For $D$ a domain in $X$, $\H(D,Y)$
denotes all holomorphic maps from $D$ to $Y$ and $C(\ol D,Y)$ denotes continuous maps from  $\ol D$ to $Y$.

%and $GL(X)$ is the group of invertible elements in $\L(X)$.$D$ to $Y$

\begin{Definition}\label{def:JBT}
    A \JBT\ is a complex Banach space $Z$ with a real trilinear
    mapping $\tp\cdot\cdot\cdot:Z\times Z\times Z \to Z$ satisfying
    \begin{enumerate}
        \item[(i)] $\tp xyz$ is complex linear and symmetric in the
        outer variables $x$ and $z$, and is complex anti-linear
        in $y$,
 \item[(ii)] The map $z\mapsto \tp xxz$, denoted $\D xx$, is
                Hermitian, $\sigma(\D xx)\ge 0$ and
        $\norm{\D xx} = \norm x^2$ for all $x\in Z$, where
        $\sigma$ denotes the spectrum,
\item[(iii)] The product satisfies the \df{Jordan triple identity}, namely
            \begin{equation}\tp ab{\tp xyz} = \tp{\tp abx}yz -
              \tp x{\tp bay}z + \tp xy{\tp abz}.\label{eq:JTI}
            \end{equation}
          \end{enumerate}
  \end{Definition}
  Throughout $Z$ will be a \JBT\ and $B$ its open unit ball.
  
The Jordan triple identity \eqref{eq:JTI} implies that $i\D xx$ is a triple
product derivation and it follows that $\exp i\D xx$ is a triple
product automorphism.  As $\D xx$ is Hermitian, $\exp i\D
xx$ is then both a triple automorphism and a surjective linear isometry.
In fact, a bijective linear map on a \JBT\ is an isometry if, and only
if, it is a triple homomorphism~\cite{Kaup_RMT}.

The triple product is continuous, 
$\norm{\{x,y,z\}} \leq \norm{x} \norm{y} \norm{z}$
\cite{FR_GN} and $\norm{\tp xxx}=\norm x^3$.  Odd powers of an element $x$ exist, with
$x^{(2n+1)}:=\{x,x^{(2n-1)},x\}, n \in \N, n\geq 1$, allowing us to define
$p(x)$, for any odd polynomial $p$, leading to an odd functional
calculus.  The triple spectrum $K_x\subset [0,\norm x]$ of the element
$x\in Z$ is defined as $\{t\in \R^+: t^2 \in \sigma(\D xx)\}$ \cite{Kaup_Spect}.  The smallest closed subtriple of $Z$
containing $x$, denoted $Z_x$, is triple isomorphic to the commutative
C*-algebra $C_0(K_x)$ via a linear map $j_x$ which takes  $x$ to the
identity function on $K_x$.  In particular, if $p$ is an odd
polynomial then
$j_x \circ p = p\circ j_x$. We refer to  $C_0(K_x)$ as the local structure of $Z$ at $x$.
\smallskip

There are natural linear maps $\D xy \in \L(Z):z\mapsto \tp xyz$,
$\q x \in \L_{\R}(Z):z\mapsto \tp xzx$, and the Bergmann operators
$B(x,y)=I-2\D xy + \q x\q y \in \L(Z)$.

\begin{Example}\label{example:C2}
  $\L(H,K)$, for complex Hilbert spaces $H$ and $K$, is a \JBT\ with
  $\tp xyz = \half(xy^*z + zy^*x)$, where $y^*$ denotes the usual
  adjoint of $y$.
\end{Example}

\subsection{Algebraic identities}

The Bergmann operators defined above play an important role in
constructing \df{quasi-inverses} in a \JBT.  Given $x, y \in Z$, we
say $(x,y)$ is a quasi-invertible pair, or that the quasi-inverse
$x^y$ exists, if the Bergmann operator $B(x,y)$ is invertible in
$\L(Z)$ and then the quasi-inverse $x^y$ is defined to be $B(x,y)\inv
(x-\q x(y))$.  We note that $B(x,y)$ is invertible if, and only if, $B(y,x)$ is invertible. Invertibility of $B(x,y)$ always holds when
$x,y \in B$ (indeed, more generally, when $\norm{\D xy}<1$) and then
$x^y=(\Id - \D xy)\inv x = \sum_ {k=0}^\infty (\D xy)^k x$.  The
quasi-inverse is a crucial component in the algebraic expression of
the biholomorphic automorphisms of the unit ball (see below).

The algebraic identities listed in \cite[Appendix A]{Loos_BSD} for
finite dimensional Jordan pairs are valid for general \JBT s.  We
recall several  for later convenience:
\begin{align*}
  B(x,y+z)&=B(x,y)\,B(x^y,z) \label{JP33}\tag{JP33} \\
  B(y+z,x)&= B(z,x^y)\; B(y,x)\label{JP34}\tag{JP34}\\
  B(x,y)\inv &= B(x^y,-y) \label{JP35}\tag{JP35}\\
    B(B(u,v)x, B(v,u)\inv y)&= B(u,v)\, B(x,y)\, B(v,u)\inv \label{JP36}\tag{JP36}\\
%\end{align*}
%and add two more which can be quickly derived:
%\begin{align*}
  x^{y+z} &= (x^y)^z  \tag{JPA1}\label{JPA1}\\
  (x+z)^y &= x^y+B(x,y)\inv z^{(y^x)} \tag{JPA2}\label{JPA2}\\
  (B(x,y)z)^y&=B(x,y)z^{B(y,x)y} \tag{JPS}\label{JPS}\\
\end{align*}

For $a\in B$, the Bergmann  operator $B(a,a)$ has positive spectrum and
a unique square root  $B(a,a)^\half$ with positive spectrum which we
will denote by $B_a$.   Via the functional calculus and local structure,
the use of the identities above can be somewhat extended.  For example,
\ref{JP36} may be used with the operator $B_a = B(a,a)^\half$ in place
of $B(u,v)$; that is $B(B_a x, B_a\inv y) = B_a B(x, y) B_a\inv$.
The local structure also allows one to make certain calculations in a
commutative setting such as
\begin{equation}
  a^a = \B a\inv a.\label{eq:local1}
\end{equation}

\subsection{Automorphisms}
The defining characteristic of the
unit ball $B$ of a \JBT\ $Z$ is its transitivity under the group of
biholomorphic mappings $\Aut B$, ensuring that $B$ is a bounded symmetric
domain.   The algebraic characterisation of all bounded symmetric domains
  \cite{Kaup_RMT} yields Definition~\ref{def:JBT} and
an explicit description (\cite[4.6]{Kaup_RMT}) of the elements of $\Aut B$.   To be precise,
every $g\in\Aut B$ can be written in the form $g=Tg_a$, where
$T$ is a surjective linear isomety of $Z$ and $g_a$ is a generalised
Möbius map, or \df{transvection}, defined on $B$ by
\[ g_a(x)= a+ \B a x^{-a}.\]
Evidently, $g_a(0)=a$ and one can show that $g_a\inv= g_{-a}$.
(This factorisation of $g=Tg_a$ is unique for if $Tg_a=Sg_b$ then
  $g_{-b}S\inv T=g_{-a}$ and applying to the origin yields $a=b$ and
  then $S=T$.) For
$a\in B$, the quasi-inverse $x^{-a}$, and hence $g_a$ and every
element of $\Aut B$, is defined and
continuous beyond the unit ball, to the ball of radius $\norm a\inv$.
The quasi-inverse map $x\mapsto x^y$ is holomorphic on its domain and
its derivative at $x_0$ is given by $B(x_0,y)\inv$ (see \cite{Loos_BSD}).
It follows then that the derivative of $g_a$ at $x_0$ is $B_a
B(x_0,-a)\inv$ and so $g_a'(0)=B_a$ while $g_a'(-a)=B_a\inv$.
Note that $Tg_a = g_{Ta}T\inv$ so we may also choose to write $g\in
\Aut B $ uniquely in the form $g=g_b S$.  In this case,
$b=g(0)$.

\subsection{Tripotents}
A tripotent is an element $e \in Z$ satisfying $\{e,e,e\}=e$ and,
since $\norm{\tp eee}=\norm e^3$,
any non-zero tripotent is a unit vector. Every non-zero tripotent $e$
induces a splitting of $Z$, as
$Z=Z_0(e) \oplus Z_{\frac{1}{2}}(e) \oplus Z_1(e)$, where $Z_\lambda(e)$ is
the $\lambda$-eigenspace of $\D ee$.  Mutually orthogonal projections
of $Z$ onto $Z_0(e),\ Z_{\frac{1}{2}}(e),$ and $Z_1(e)$ are given by
$P_0(e)=B(e,e), P_{\frac{1}{2}}(e)=2(\D ee -\q e\q e)$ and
$ P_1(e)=\q e\q e$, respectively .

A tripotent is called \df{maximal} if $Z_0(e)=\{0\}$ or, equivalently,
if $B(e,e)=0$ and is called \df{unitary} if $Z_1(e)=Z$, that is, $P_1(e)=\q
e\q e=\Id$.   For any \JBT, the
set of real extreme points of $\overline B$, the set of complex
extreme points of $\ol B$ and the set of maximal tripotents all
coincide.  For details see \cite{Loos_BSD}.

Let $\Gamma$ be the set of maximal tripotents and $\Gamma_1$ be the
set of unitary tripotents of the \JBT\ $Z$.   If  $B(a,a)=0$ or
$\q a \q a=\Id$ then
consideration of $0=B(a,a)a$, or $a=\q a^2 a$  in the local triple $Z_a$ shows $a$ to be
a tripotent and so:
\begin{Lemma}\label{extreme}
  $\Gamma =\{a\in \ol B : B(a,a)=0\}$ and $\Gamma_1=\{a \in\ol B: \q
  a\q a= \Id\}$.
\end{Lemma}

An important  consequence of Lemma \ref{extreme} is that both $\Gamma$ and $\Gamma_1$ are closed.
We say $x,y \in Z$ are orthogonal,
$x \bot \ y$, if $x \square y = 0$ (equivalently $y \square x =
0$). In particular, if $c$ and $e$ are orthogonal tripotents then
$c+e$ is also a tripotent. This gives a partial ordering on the set,
$M$, of all tripotents in $Z$ as follows.

\begin{Definition}
  For tripotents $c$ and $e$ we say $c<e$ if $e-c \in M$ and
  $(e-c) \bot \ c$.
\end{Definition}

Maximality of a tripotent with respect to this ordering is consistent with
the notion of maximal tripotent given previously in the sense that a
maximal tripotent is order maximal.  A tripotent $e$ is
minimal if $Z_1(e)=\C e$. $Z$ is said to have finite rank $r$ if every
element $z \in Z$ is contained in a subtriple of (complex) dimension
$\leq r$, and $r$ is minimal with this property. We say $x$ is a rank $k$ element if $Z_x$ is $k$-dimensional. The rank one
\JBT s are the Hilbert spaces.  Other finite rank examples are
sub-triples of $\L(H,K)$ where either $H$ or $K$ is finite
dimensional, and spin factors.  If $Z$ has finite rank $r$, a frame is
a set $\{e_1, \ldots, e_r\}$ of non-zero pairwise orthogonal minimal
tripotents and every $z \in Z$ has a unique spectral decomposition,
called its Peirce decomposition, as
$z=\lambda_1 e_1+\cdots +\lambda_r e_r$, for some frame
$\{e_1, \ldots, e_r\}$ and scalars
$0\leq \lambda_1 \leq \ldots \leq \lambda_r=\|z\|$.  See \cite{Din_Sch,Kaup_I} for details.

A \JBT\ may not have any tripotents but if it is a dual Banach space
(for example, the bidual of a \JBT\ is a \JBT\ \cite{Dineen_CHVF}) then
the Krein-Millman theorem implies the existence of maximal
tripotents.   Finite rank \JBT s have maximal tripotents and indeed
are reflexive Banach spaces.

\section{Invariance under $g\in \Aut B$}

As mentioned above, biholomorphic automorphisms of the open unit ball extend
continuously to a neighbourhood of $\ol B$.    These maps do not, however, preserve the set of
tripotents generally.  

\begin{Example}\label{tripotent} The commutative C*-algebra, $\C^2$, is a 
\JBT\  via the product $\tp fgh = f\ol g h$ (corresponding to the maximum norm).  Take $e$ to be the tripotent $(1,0)$ and
$a=(\half,\half)$ and notice that for $g_a \in \Aut B$, $g_a(e)=(1,\half)$
is not a tripotent.
\end{Example}

Nonetheless, the automorphisms of $B$ do act invariantly
both on the set of maximal tripotents and on the set of unitary
tripotents (when these are non-empty).  In order to prove this we require the
following results establishing key identities involving the Bergmann
operators.

\begin{Proposition}\label{Bergmann_1}
  Let $a\in B$ and $b\in\ol B$.  Then \[ B(g_a(b),g_a(b))= B_a B(b,-a)\inv B(b,b)
    B(-a, b)\inv B_a.\]
\end{Proposition}

\begin{Proof}
Recalling that $g_a(b)= a+ B_ab^{-a}$ we proceed to expand as follows.
\begin{align*}
  B(\g ab, \g ab) =& B(a+ \B ab^{-a}, a+ \B ab^{-a}) \\
     \overset{\eqref{JP34}}=& \underlabel{R}{B(\B ab^{-a} , (a+ \B a(b^{-a}))^a)} \, . \,
                            \underlabel{S}{B(a, a+  \B a(b^{-a}))}  \\
\end{align*}
Focusing on the Bergmann operator $R$,  we have
\begin{align*} (a+\B a(b^{-a}))^a \overset{\eqref{JPA2}}=& a^a+ B(a,a)\inv (\B a
                                                          (b^{-a}))^{(a^a)}\\
                          \overset{\eqref{eq:local1}}=\ & \B a\inv[a+ \B
a\inv(\B a (b^{-a}))^{\B a\inv a}] \\ \overset{\eqref{JPS}}=& \B a\inv[a+ \B a\inv
                                                      \B a (b^{-a})^a]= \B a\inv (a+b)
\end{align*}
and consequently
\begin{align*}
  R=& B(\B a b^{-a}, \B a\inv(a+b)) \\
                          \overset{\eqref{JP36}}=&\B a\, B(b^{-a},
                                                   a+b) \, \B a\inv
  \\
  \overset{\eqref{JP33}}=& \B a \, B(b^{-a},a)\, B(b,b)\,  \B a\inv \\
  \overset{\eqref{JP35}}=& \B a\,
                           B(b,-a)\inv \, B(b,b)\, \B a \inv.
\end{align*}
A similar expansion gives
\[ S = \B a B(-a, b)\inv \B a\]
and we are done.
\end{Proof}

\begin{Lemma}\label{lem:isom}
  Let $a, b \in B$.  Then
  \begin{enumerate}[(i)]
  \item $k(a,b)=  \B {g_a(b)}\inv \B a  B(b,-a)\inv \B  b$ is a
    surjective linear
    isometry of $Z$,
    \item $g_a g_b= g_{\g ab} k(a,b)$.
  \end{enumerate}
\end{Lemma}

\begin{Proof}
  Consider the biholomorphic map on $B$ given by $k(a,b)=g_{-\g ab} g_a
  g_b$. Since this map fixes $0$, the Schwarz Lemma guarantees it is a
  linear isometry which agrees with its derivative at the origin and
  this, 
  via the chain rule, is given by \[g_{-\g ab}'(\g
  ab).g_a'(b).g_b'(0)= \B {\g ab}\inv \B aB(b,-a)\inv \B b.\]
\end{Proof}

We now provide proofs of the invariance of the set of maximal tripotents
$\Gamma$, and the set of unitary tripotents $\Gamma_1$ under
automorphisms.  The result for $\Gamma$ is effectively provided in
\cite[2.4(i)]{MR2834260} using the characterisation of maximal
tripotents as complex extreme points.  Here we present an algebraic
proof which follows immediately from Proposition~\ref{Bergmann_1}. 

\begin{Theorem}\label{thm:GammaInv}
  For $g\in \Aut B$, $g(\Gamma)=\Gamma$, \  if $\Gamma\neq \emptyset.$
\end{Theorem}

\begin{Proof}
  Write $g=Tg_a$ where $T$ is a surjective linear isometry and $g_a$
  is a Möbius map.  Surjective linear isometries (i.e. triple product automorphisms) not
  only preserve tripotents, but also preserve the maximality of a
  tripotent.   Indeed, $f\in\ol B$  is a maximal tripotent precisely
  when $B(f,f)=0$ and, as $T$ is a triple automorphism,  this
  coincides with $T B(f,f) T\inv = B(Tf, Tf)$ vanishing and $Tf$ being
  a maximal tripotent.

  Thus we must show that for $a\in B$ and a
  maximal tripotent $e$, we have $g_a(e)$ is a maximal tripotent.
  Since $e$ is a maximal tripotent, $B(e,e)=0$ and from
  Proposition~\ref{Bergmann_1} we have $B(g_a(e),g_a(e))=0$ and thus
  $g_a(e)$ is a maximal tripotent.   This proves inclusion, while
  equality follows by invertibility and $g_a\inv = g_{-a}$.
\end{Proof}

\begin{Proposition}\label{prop:swap}
  Let $a, b \in B$.  Then $\g ab = \tilde k \g ba$, where $\tilde k$ is
  a surjective
  linear isometry.
\end{Proposition}

\begin{Proof}
  By inversion of Lemma~\ref{lem:isom}(ii), we can say \[g_{-b}g_{-a}=
  k(a,b)\inv g_{-\g ab} = k(a,b)\inv g_{g_{-a}(-b)} \] and swapping the
  roles of $a$ and $-b$, we have
  \[ g_a g_b = k(-b, -a)\inv g_{\g ba}.\]
  Comparison with Lemma~\ref{lem:isom}(ii) yields $g_{\g ab} k(a,b)=
  k(b,a)\inv g_{\g ba}$.  Apply this to the origin to gain the result,
  together with $\tilde k = k(b,a)\inv$.
\end{Proof}

\begin{Theorem}\label{thm:Gamma1Inv}
  For $g= \in \Aut B$,  $g(\Gamma_1)=\Gamma_1$, \   if $\Gamma_1\neq \emptyset$.
\end{Theorem}

\begin{Proof}
  Again write $g=Tg_a$ where $T$ is a surjective linear isometry and $g_a$
  is a Möbius map.   As in the proof for maximal tripotents, invariance with respect to
  the linear part $T$ is immediate and we need only show that the
  image of a unitary tripotent under the transvection $g_a$
  is a unitary tripotent.   This equates to proving 
  $\q{g_a(u)}\q{g_a(u)}=\Id$ when $\q u\q u=\Id$.

  Let $t\in(0,1)$ so that $tu \in B$.  
  Proposition~\ref{prop:swap} allows us to  write $g_a(tu)= k_t g_{tu}(a)$,
  for some surjective linear isometry $k_t$.  As $k_t$ is a triple
  automorphism it follows that $\q{g_a(tu)} =\q{k_tg_{tu}(a)}
  =k_t\q{g_{tu}(a)}k_t\inv$ and so
  \[\Id - \q{g_a(tu)}\q{g_a(tu)} = k_t(\Id
    -\q{g_{tu}(a)}\q{g_{tu}(a)})k_t\inv.\]
  In particular, \begin{equation}\norm{\Id - \q{g_a(tu)}\q{g_a(tu)}}=\norm{\Id
    -\q{g_{tu}(a)}\q{g_{tu}(a)}}.\label{eq:norm}
\end{equation}
Being unitary, the tripotent $u$ is also maximal and the Bergman
operator $B(u,u)=0$.  Triple product continuity implies then that
$B(tu,tu)\to 0$.  This convergence passes to the square root.  Indeed,
 by~\cite[Lemma 3.4]{Kaup_Sauter},
 $\norm{B_{tu}}=\norm{B_u-B_{tu}}\le 2\sqrt{1-t^2}$ and so $B_{tu}\to
   0$ as $t\to 1$.  Since $a^{-tu} \to a^{-u}\in Z$, we can say that
   $g_{tu}(a)= tu+ B_{tu}a^{-tu} \to u$ as $t\to 1$.  Again from continuity of
   the triple product, $\q{g_{tu}(a)} \to \q u$ and hence $\norm{\Id
    -\q{g_{tu}(a)}\q{g_{tu}(a)}} \to 0$.   Combine this with
  \eqref{eq:norm} and the convergence of $g_a(tu)$ to $g_a(u)$ to
  conclude that $\q{g_a(u)}\q{g_a(u)} =\Id$ as required.
\end{Proof}

 \section{Russo-Dye type results for JB*-triples}

 The classical Russo-Dye Theorem states that the closed unit ball of
 a unital C*-algebra is the closed convex hull of its extreme points \cite{MR0193530}.  
Variations of the second statement in Theorem~\ref{triple_russo-dye} below appear in \cite{Mackey_homotopes} and 
 \cite{MR2376600}, but here we do not require that $Z$ is a JBW*-triple
 (i.e. has a predual).  We also provide a unified proof of both statements below, based on a
 technique used by Harris \cite{Harris_bshd}.
 \begin{Theorem}\label{triple_russo-dye}
   Let $Z$ be a JB*-triple with open unit ball $B$.
   \begin{enumerate}
   \item If $Z$ contains a maximal tripotent then $\ol B =\cch (\Gamma)$.
     \item If $Z$ contains a unitary tripotent then $\ol B =\cch (\Gamma_1)$
   \end{enumerate}
 \end{Theorem}

% Yes  [[ Is this ok?  Reference to Harris??]]

 \begin{Proof}
   Let $b \in B$.   For (1), choose $a \in \Gamma$ and, respectively for (2), $a\in \Gamma_1$.
   Define  $h(\lambda)= g_b(\lambda a)$ to gain a $Z$-valued
   holomorphic function on the disc of radius $\frac 1{\norm b} >1$ in
   $\C$.  The mean value property for holomorphic functions implies
   \begin{align*}
   b = h(0) &= \frac 1{2\pi} \int_0^{2\pi} h(e^{i\theta}) d\theta
     \\
            &=\frac 1{2\pi} \int_0^{2\pi} g_b(e^{i\theta }a) d\theta.
   \end{align*}
   For each
   $\theta\in\R$, $e^{i\theta} a \in \Gamma$ (resp. $\Gamma_1$) and, by Theorem~\ref{thm:GammaInv}
   (resp. Theorem~\ref{thm:Gamma1Inv}), so is $g_b(e^{i\theta}a)$.   Thus $b \in \cch (\Gamma)$
   (resp. $\cch(\Gamma_1)$) as stated.  
 \end{Proof}

 In fact, more general Russo-Dye extensions are possible.  For any $v\in
 \del B$, we let $G_v := \Aut B (v)$ be the orbit of $v$ under the
 automorphisms of $B$ and we let $G_v^0 := \Aut B^0 (v)$ be the orbit
 of $v$
 under the automorphisms in $\Aut B^0 $, 
 the connected component of the identity in $\Aut B$.
Note that every Möbius map $g_a$ and every unimodular rotation
 $e^{i\theta}$ lies in $\Aut B^0$. Replacing $a$ in the above proof with $v \in \del B$, means  $g_b(e^{i\theta }v) \in G^0_v,$ for all $b \in B$ and $\theta \in \R$ and yields the following.
 \begin{Corollary}\label{general_hull}
   Let $Z$ be a \JBT\ with open unit ball $B$.  For any $v\in\del B$,
   $$\ol B = \cch (G^0_v)= \cch (G_v).$$
 \end{Corollary}
 
 We will see later in Corollary~\ref{last} that if $Z$ is a finite rank triple then 
 $\Gamma \subseteq  \ol {G^0_v},$ for all $v \in \del B$.
 
 \section{Boundary subsets that are determining for holomorphic functions}
 
An immediate Corollary of Theorem \ref{triple_russo-dye} is the result below (and equation  \eqref{eq:value} in particular) showing that holomorphic functions on $B$ having a continuous extension to $\ol B$
are already determined by their values on  $\Gamma$ or
$\Gamma_1$, when these are non-empty. This can be described by saying that the sets  $\Gamma$
and $\Gamma_1$ are {\it determining} for holomorphic functions on $
B$. In fact, there are many such determining sets in the boundary of a bounded symmetric domain but for clarity we begin with $\Gamma$ and $\Gamma_1$. 
 
  \begin{Proposition}\label{determining}
   Let $Z$ be a \JBT\ with open unit ball $B$ and $X$ be any Banach space. Let $f:B \mapsto X$ be holomorphic on $B$ with a continuous extension to $\ol B$. 
 \smallskip
   
   \begin{enumerate}
   \item If $Z$ contains a maximal tripotent then $ f(\ol B) \subseteq\cch (f(\Gamma))$
   and $$\sup\{\|f(z)\|:z \in \ol B\}=\sup\{\|f(z)\|:z \in \Gamma\}.$$
     \item If $Z$ contains a unitary tripotent then $f(\ol B) \subseteq\cch (f(\Gamma_1))$ and $$\sup\{\|f(z)\|:z \in \ol B\}=\sup\{\|f(z)\|:z \in \Gamma_1\}.$$
      \end{enumerate}  
   \end{Proposition}
 \begin{Proof} The two versions have similar proof, so we only present
   that of (1).  Suppose then
   that $\Gamma \neq \emptyset$ and then let $a
   \in \Gamma$.  Let $b \in B$ be arbitrary.  Applying the proof of Theorem \ref{triple_russo-dye}  to   
the map  $h(\lambda)= f(g_b(\lambda a))$ and using the  mean value property gives 
   \begin{align}
  f(b) =\frac 1{2\pi} \int_0^{2\pi} f(g_b(e^{i\theta }a)) d\theta. \label{eq:value}
   \end{align}
   It follows that $$\|f(b)\|\leq \frac 1{2\pi} \int_0^{2\pi}\sup\{\|f(z)\|:z \in \Gamma\} d\theta =\sup\{\|f(z)\|:z \in \Gamma\}$$ and we are done.  
   \end{Proof}
\begin{Remarks}We note that, unlike in finite dimensions, $\sup\{\|f(z)\|:z \in \ol B\} $ may not actually be achieved. However, the above means that if $f$ is bounded on $\Gamma$ or
$\Gamma_1$, it must also be bounded on $\ol B$ and thus, if $f$ is
unbounded then all suprema above must be infinite.
\end{Remarks}
 Just as Theorem~\ref{triple_russo-dye} leads to
 Proposition~\ref{determining}, we have the following
 from Corollary~\ref{general_hull}.
 
 \begin{Proposition}\label{vorbits}
      Let $Z$ be a JB*-triple with open unit ball $B$ and $X$ be any Banach space. Let $f:B \mapsto X$ be holomorphic on $B$ with a continuous extension to $\del B$.  For any $v\in \del B$ we have
    \[ f(\ol B) \subseteq \cch f(G_v^0)\]
    and
    \[\sup \{\norm{f(z)}:z\in \ol B\} =   \sup\{ \norm{f(z)}: z\in G_v^0\}.\]
    
    In particular, if $f$ is unbounded then both suprema are infinite.
 \end{Proposition}
 
 We recall that if a \JBT\ has finite rank, then the boundary of its unit ball is the disjoint union of holomorphic boundary components defined as follows.

\begin{Definition}\cite[4.1]{Kaup_Sauter}\label{bdycomp}
  A non-empty set $A \subset \overline{B}$ is a holomorphic boundary
  component of $B$ if $A$ is minimal with respect to the fact that,
  for all $f \in \mathcal{F}=\{f: \Delta \rightarrow Z$ holomorphic
  with $f(\Delta) \subset \overline{B}\}$, either
 $$f(\Delta) \subset A\ \hbox{or}\ f(\Delta) \subset \overline{B} \setminus A.$$
\end{Definition}
By
replacing $\mathcal{F}$ in the above definition with the set of all
complex affine maps $:\Delta \rightarrow \overline{B}$ we get the
definition of (complex) affine boundary components.

\begin{Remarks}\label{domain}It follows that if $D$ is a domain and
  $g:D\rightarrow \overline B$ is a holomorphic map then $g(D)$ must
  lie in a single such boundary component.
\end{Remarks}

The following shows that holomorphic and affine boundary components
coincide in the finite rank case and each is determined by a unique tripotent.

\begin{Theorem}{\normalfont \cite[4.2, 4.3, 4.4 ]{Kaup_Sauter}}\label{Loos-Boundaries}
  Let $Z$ be a finite rank \JBT\ with open unit ball $B$. The
  following hold.
\begin{enumerate}
\item[(i)] Holomorphic and affine boundary components coincide and are
  precisely the sets $$K_e=e+B_0 (e)$$ where $e$ is a tripotent and
  $B_0 (e)=B \cap Z_0 (e)=P_0(e)(B)$.
\item[(ii)] The map $e \rightarrow K_e$ is a bijection between the
  set, $M$, of tripotents in $Z$ and the set of boundary components of
  $B$, with $x \in K_e$ if, and only if,
  $e=\lim_{n\rightarrow \infty}x^{(2n+1)}.$
\item[(iii)] $\ol K_e=\bigcup_{d\geq e}K_d.$  
\end{enumerate}
\end{Theorem}
Henceforth we refer simply to boundary component and write $K_x$ for the boundary component of $x$. We
note that $x$ is an extreme point if, and only if, $K_x=\{x\}.$
Indeed, if a boundary component $K_x$ contains an extreme point $v$
then $K_x=\{v\}.$

\begin{Definition}\label{rank}We define the rank of a 
  boundary component $K_e$ in a finite rank triple $Z$ as the rank of
  the JB*-subtriple $Z_0(e)$ (with rank zero if $Z_0(e)=\{0\}$).
\end{Definition}

Since $K_e=e+B_0 (e)$, where $B_0 (e)=B \cap Z_0 (e)$, then $K_e$ is
biholomorphically equivalent to the bounded
symmetric domain $B_0 (e)$, which is the open unit ball of the \JBT\
$Z_0 (e)$. In other words, the rank of the boundary component $K_e$ in the sense of Definition~(\ref{rank}) is prescisely its rank as a bounded symmetric
domain. Moreover, if a tripotent $e$ is a rank $k$ element ($Z_e$ is $k$ dimensional), then its
boundary component $K_e$ is rank $n-k.$ This means that if $Z$ is rank
$n$ then it has boundary components of rank
$ k \in \{0,1,\ldots, n\}.$ The only rank $n$ boundary component is
$B (=K_0)$ itself, the rank zero components are singletons given by
the maximal extreme points for which $B_0(e)=\{0\}$, and there are
components of all rank $k \in \{0,1,\ldots, n\}.$

Theorem~\ref{thm:GammaInv}, which underpins the Russo-Dye extensions
in Theorem~\ref{triple_russo-dye} and Corollary~\ref{general_hull},
showed that the extreme points, or rank zero boundary components, are
invariant under elements of $\Aut B$. In fact,
Theorem~\ref{thm:GammaInv} is a special case of a more general result,
namely, that every automorphism of $B$ maps rank $k$ boundary
components onto rank $k$ boundary components, for all
$k \in \{0,1,\ldots, n\}.$ We note that this holds true despite the fact that
automorphisms do not generally preserve the rank of individual
elements, nor must automorphisms even map tripotents to tripotents, as
seen in Example~\ref{tripotent} above.

For $g \in \Aut B$, we recall that $g:\del B\mapsto\del B$ and since
$g$ extends to a holomorphic map on an open neighbourhood of $\ol B$,
it follows that $g|_{K_v}$ is holomorphic on $K_v$.

\begin{Proposition}\label{preserve-rank}
  Let $Z$ be a finite rank \JBT. Let $v \in \ol B$ and $g \in \Aut B$. Then
\smallskip

\begin{enumerate}
\item $g(K_v)=K_{g(v)};$\\
\item  $\rank(K_v)= \rank(K_{g(v)})$.
\end{enumerate} 
  \end{Proposition}

  \begin{Proof}
    Let $v \in \ol B$ be arbitrary and $K_v$ be its boundary
    component.  \newline Case (i): Suppose $v$ is not extreme, so
    $K_v$ is biholomorphically equivalent to a non-trivial domain. Then
    $g(K_v) \cap K_{g(v)}\neq \emptyset$ (as it contains $g(v)$), so
    from Remarks~\ref{domain} $g(K_v) \subseteq K_{g(v)}.$ Then
    $K_{g(v)}$ is not a singleton, so $g(v)$ is not extreme, and the
    above argument applied to $g(v)$ and $g^{-1}$ gives
    $g^{-1}(K_{g(v)})\subseteq K_{v}$. Thus $g(K_v)= K_{g(v)}.$
    \newline Case (ii): Let $v$ be extreme. Then $g(v)\in K_u$, for
    some tripotent $u$ and if $u$ is not extreme, applying (i) to $u$
    and $g^{-1}$ gives $g^{-1}(K_u)=K_{g^{-1}(u)}.$ Then
    $v \in K_{g^{-1}(u)}$ is extreme, so
    $\{v\}=K_{g^{-1}(u)}=g^{-1}(K_u)$. This is impossible as $g^{-1}$
    is injective and $K_u$ is a non-trivial domain.  In other words,
    $v$ extreme implies $u$ extreme and hence $g(v)$ is extreme,
    completing (1) above.

To prove (2). Let $v \in \ol B$. From (ii), $v$ is extreme if, and only if, $g(v)$ is extreme, in which case, $\rank(K_v)= \rank(K_{g(v)})=0.$ We assume therefore that $v, g(v)$ are not extreme. There exist tripotents $e,f$ with (and from (1))
$$K_v=K_e=e+B_0(e)\ \hbox{and}\ g(K_v)=K_{g(v)}=K_f=f+B_0(f),$$
where  $B_0(e)$ and $B_0(f)$ are the open units balls of (non-zero) $JB^*$-triples  $Z_0 (e)$ and $Z_0 (f)$ (respectively). We define $$h: B_0 (e)\mapsto  B_0 (f)\ \hbox{given by}\ 
h(z)=g(z+e)-f.$$ Clearly, $h$ is holomorphic with holomorphic inverse
$$h^{-1}: B_0 (f)\mapsto  B_0 (e)\ \hbox{given by}\ 
h^{-1}(w)=g^{-1}(w+f)-e.$$
In other words, the open unit balls of $Z_0 (e)$ and $Z_0 (f)$ are biholomorphically equivalent (under $h$). By \cite{KaupUpmeier_biholo}, two Banach spaces are linearly isometric if, and only if, their open unit balls are biholomorphically equivalent, so $Z_0 (e)$ and $Z_0 (f)$ are linearly isometric. A linear isometry of triples preserves the triple rank, so $Z_0 (e)$ and $Z_0 (f)$ have the same rank as triples and hence $K_v=K_e$ and  $K_{g(v)}=K_f$ have the same rank as boundary components. 
\end{Proof}
\smallskip

\begin{Corollary}
Let $Z$ be a \JBT\ of finite rank $n$. For each $k \in \{0,\ldots, n-1\}$ there is a determining set in $\del B$ whose points all lie in rank $k$ boundary components.
\end{Corollary}

\begin{Proof} Fix $0\leq k\leq n-1$. Choose any tripotent $e$ of rank
  $n-k$. The boundary component $K_e$ has rank $k$. Pick any
  $v \in K_e$. Proposition~\ref{preserve-rank} above proves that for
  all $g \in G=\Aut B$, $\rank(K_{g(v)})=\rank(K_v)$. In other words, each element in
  $G_v$ (and hence in $G_v^0$) lies in a rank $k$ boundary component.
\end{Proof}

We note this does not mean that $G_v^0$ or $G_v$ (for $v \in \del B$)
lies in any one boundary component. For example, if $k=0$, the rank
$k$ boundary components are the extreme points in $\Gamma$, and this
set is generally not connected in the infinite dimensional case.
Also of relevance to later results is the fact that the
  determining sets $G_v$ above are not closed in general.  To
  illustrate, let us return to Example~\ref{example:C2}.

  \begin{Example}
    Take $v=(1,\half)$ in the boundary of the unit ball of $Z=\C^2$ for the maximum norm.  Then
    $G_v^0=\T\times \Delta$ while $G_v=(\T \times \Delta) \cup
    (\Delta\times \T)$.  The
    holomorphic boundary component of $v$ is $K_v=\{1\}\times\Delta$, and this
    equals $K_e$ where $e=(1,0)$ is the unique tripotent in $K_v$.  
    The closure of $G_v^0$  is a proper subset of the boundary which
    properly contains the set of extreme points $\Gamma=\T\times\T$.
    Each extreme point is not only a maximal tripotent, but also
    unitary.   These last two facts will be reflected in
    Corollaries~\ref{maximalisunitary} and \ref{last}. 
  \end{Example}

We have shown that there are many determining sets in the boundary of
a bounded symmetric domain, for example,
$\Gamma_1, \Gamma, G_v, G_v^0$, for arbitrary $v \in \del
B$. Nonetheless, we will show that the role played by the set, $\Gamma$, of extreme
points remains special.  We recall that in finite dimensions, the
Bergmann-Shilov boundary of $B$ is defined as the minimal closed
subset of $\ol B$ on which every $f:\ol B\mapsto \C$ which is
holomorphic on $B$ and continuous on $\ol B$ achieves its maximum
modulus. Moreover, the Bergmann-Shilov boundary in finite dimensions
is exactly the set of extreme points, $\Gamma$ \cite[Theorem
6.5]{Loos_BSD}.

While Propositions~\ref{determining} and \ref{vorbits} can already be
viewed as partial extensions of Bergmann-Shilov type behaviour to
infinite dimensions, in the next section we will prove that for finite
rank triples, the set of extreme points $\Gamma$ is the key determining set and is exactly a
Bergmann-Shilov boundary relative to such holomorphic functions.
   
\section{An infinite dimensional analogue of the Bergmann-Shilov
  boundary} 

Let $Z$ be a finite rank \JBT.  We recall now an equivalent norm
defined on $Z$ by means of the spectral decomposition \cite[Section
9.2]{Din_Sch}.  Namely, for $x \in Z$ we have
$x=\alpha_1 e_1+\cdots +\alpha_n e_n,$ with
$\|x\|=\alpha_1\geq\alpha_2\cdots \geq \alpha_n\geq 0,$ and
$e_1, \cdots, e_n$ a frame of mutually orthogonal minimal
tripotents. To each minimal tripotent, \cite{Friedman-Russo_predual},
$e_i$ there exists a unique $\phi_i \in Z_*$ (the unique predual of
$Z$) such that $\phi_i$ is extreme in $\ol B_{Z_*} $ and
$\phi_i(e_i)=1$.  This allows us to define an inner product, often called the
\df{algebraic inner product}, on $Z\times Z$ by
$\ip xy_a:=\sum_{i=1}^n \alpha_i \ol{\phi_i(y)}$ and
$\|x\|_a^2:=\ip xx_a=\sum_{i=1}^n \alpha_i^2$
and \begin{align}\label{algnorm}\|x\|^2\leq \|x\|_a^2\leq
  n\|x\|^2.\end{align}

We quote the following version of the Maximum Modulus Principle
\cite[Cor. 2.2]{MR2022955}.

\begin{Theorem}[Maximum Modulus]\label{maxmod}
  Let $D$ be a domain in a complex Banach space $Z$ and
  $f\in \H(D,X)$, where $X$ is a Banach space.  Then
  $f$ satisfies the maximum modulus principle, namely, if
  $\norm{f(z)}$ achieves a maximum at any point of $D$ then $\norm{f(z)}$ is
  constant.
\end{Theorem}
      
\begin{Proposition}\label{prop:A}
  Let $Z$ be a finite rank \JBT, $X$ be a Banach space
  and $f$ be a holomorphic map on a neighbourhood of $\ol B$ into $X$.  If
  $\norm{f(z)}$ achieves a maximum on $\ol B$ then it must achieve
  this maximum on $\Gamma$.
\end{Proposition}

\begin{Proof}
  Suppose $\norm{f(z)}$ achieves a maximum $M$ at $w\in\ol B$.  If $w\in
  B$ then by the maximum modulus principle above, $\norm{f(z)}$ is constant on
  $B$ and hence on $\ol B$, thereby achieving its maximum at every
  element of $\Gamma$.   Thus we may assume $w\in \del B$.   If $w$ is an
  extreme point we are done, so we may assume that $w$ is an interior
  point of its boundary component $K_w$.

  Since $Z$ is finite rank, there is a unique tripotent $e$ such that
  $K_w=K_e = e+ (Z_0(e)\cap B)$.  We define a holomorphic map $h$ on
  $D:=Z_0(e)\cap B$ by $h(z)=f(e+z)$.  Clearly, $\|h(z)\|$ is bounded by $M$
  and achieves this bound at $z_0\in D$ where $w=e+z_0$.  Again
  by Theorem~\ref{maxmod}, $\norm{h(z)}$ is constant on $D$ and so $\norm{f(z)}$ is constant on
  $K_e$ and $\|f(e)\|=\|f(w)\|=M$. By continuity, $\norm{f(z)}$ is then constant
  on $\ol{K_e}$. As $\ol{K_e}= \cup_{d \ge
    e} K_d$, $\norm{f(z)}$ is also constant on $K_d$, where $d$ is any
  tripotent that majorises $e$.  Using the spectral decomposition, we
 can construct a maximal tripotent $e'$ that majorises $e$. By maximality,  $K_{e'}=\{e'\}$ so we have $e'\in \Gamma$ with $\|f(e')\|=\|f(e)\|=M$ as required.
\end{Proof}

The hypothesis in Proposition~\ref{prop:A} that $f$ be holomorphic on
a open neighbourhood of $\ol B$ was chosen to ensure that the mapping
$h(z)=f(e+z)$ defined on $D=Z_0(e)\cap B$ is itself holomorphic; where
$K_e= e+ (Z_0(e)\cap B)$ is the boundary component on which $f$
achieves maximum norm. In finite dimensions, it suffices for
$f:B\mapsto \C$ that $f\in \H(B,\C)\cap C(\ol B,\C)$.  In that case, we use
compactness of $\ol B$ and uniform convergence of the maps
$h_n(z)=f((\frac{n-1}{n})e+z)$, which are clearly now holomorphic on
$B$, to argue that $h$ is holomorphic as it is a uniform limit of
holomorphic functions. Of course, in infinite dimensions $\ol B$ is
not compact so we need other tools. These tools consist of a different
topology on the space of all maps $\H(B,\C)$, known as the
compact-open topology and denoted here by $\tau$, together with a
suitable $\tau$ analogue of Montel's theorem and knowing that $\tau$
is complete. We therefore use the following results, cf.
\cite[Proposition 2.4 and Theorem 2.10]{MR2022955}.

\begin{Theorem}\label{complete} Let $X,Y$ be Banach spaces and $D$ be a domain in $X$. The space $\H(D,Y)$ is complete with respect to $\tau$.
\end{Theorem}
We write $\|f\|_D:=\sup \{\norm{f(z)}:z\in D\}.$

\begin{Theorem}\label{montel} Let $X,Y$ be Banach spaces and $D$ be a domain in $X$.
Then the set ${\mathcal F}_M=\{f \in \H(D,Y):\|f\|_D\leq M< \infty\}$ is relatively compact with respect to $\tau$ if, and only if, each orbit 
${\mathcal F}_M(x)$  is relatively compact in $Y$, for $x\in D$.
\end{Theorem}

\begin{Theorem}\label{hol-on-B}
  Let $Z$ be a finite rank \JBT\  with ball $B$. 
  Let $f:B\mapsto \C$ be holomorphic on $B$ with a continuous extension to $\del B$.  If
  $|f(z)|$ achieves a maximum on $\ol B$ then it must achieve
  this maximum modulus on $\Gamma$.
\end{Theorem}

\begin{Proof}
 Assume $\absval{f(z)}$ achieves a maximum $M$ at $w\in\ol B$ and repeat the proof of Proposition~\ref{prop:A}. The only part of that proof that requires adapting to $f\in \H(B,\C)\cap C(\ol B,\C)$ is to prove that $h:D=Z_0(e)\cap B\mapsto \C$ given by
 $h(z)=f(e+z)$ is holomorphic, where $w=e+z_0$, as before.
 For $n \in \N$, define the holomorphic map $h_n:D \mapsto \C$ by
  $h_n(z)=f((\frac{n-1}{n})e+z)$. Using terminology from Theorem~\ref{montel} with $Y=\C$,  ${\mathcal F}_M(x)
  \subset B(0,M) \subset \C$  is bounded and is hence relatively compact in $\C$, for each $x\in D$. Theorem~\ref{montel} therefore implies that ${\mathcal F}_M\subseteq \H(D,\C)$ is  $\tau$-relatively compact in $\H(D,\C)$. The sequence $(h_n)_n$ in ${\mathcal F}_M$ must therefore have a $\tau$-convergent subnet (indexed by $\alpha$, say) $(h_{n_\alpha})_\alpha$, converging to a $\tau$-limit $k.$ Theorem~\ref{complete} then implies that $k \in \H(D,\C)$. 
 Since $\tau$-converence implies pointwise convergence, we have in particular that, for all $z\in D$,
 $$k(z)=\lim_\alpha h_{n_\alpha}(z)=\lim_\alpha f((\frac{n_\alpha-1}{n_\alpha})e+z).$$ Continuity of $f$ to $\del B$ then gives $k(z)=f(e+z)=h(z)$. In other words, $h=k$ and therefore $h$ is holomorphic on $D$.  The rest of the proof then continues exactly as in Proposition~\ref{prop:A}.
 \end{Proof}

It follows from Proposition~\ref{prop:A}, Theorem~\ref{montel} and the proof of Theorem~\ref{hol-on-B}  that $\C$ in Theorem~\ref{hol-on-B} can be replaced by any finite dimensional  Banach space.  
\begin{Corollary}
  Let $Z$ be a finite rank \JBT\  with ball $B$ and $X$ be any finite dimensional  Banach space. 
  Let $f:B\mapsto X$ be holomorphic on $B$ with a continuous extension to $\del B$.  If
  $\norm{f(z)}$ achieves a maximum on $\ol B$ then it must achieve
  this maximum  on $\Gamma$.
\end{Corollary}

The following result is adapted from \cite[Theorem 6.5]{Loos_BSD}.

\begin{Proposition}\label{prop:B}
  Let $Z$ be a finite rank \JBT\ and $X$ be a Banach
  space. For $e\in\Gamma$ there exists a holomorphic map $h:Z\to X$
  such that $h$ achieves its maximum norm on $\ol B$ only at $e$.
\end{Proposition}

\begin{Proof}
  Fix $e\in \Gamma$ and $v\in X$ with $\norm v=1$. Let $\ip
  \cdot\cdot_a$ denote the algebraic inner product on $Z$ defined
  above.  Define a holomorphic function $h:Z\to X$ by $h(z)=
  \half\left( 1+ \frac{\ip ze_a}{\ip ee_a}\right) v$. For $z\in \ol
  B$,
  \begin{equation} \norm{h(z)} = \half \absval{1+ \frac{\ip ze_a}{\ip ee_a}} \le
    \half \left(  1+ \frac{\absval{\ip ze_a}}{\ip ee_a} \right).\label{eq:1}\end{equation}
    
    The Cauchy-Schwarz inequality gives 
      \begin{equation}\label{CS}
      \absval{\ip ze_a} \le \norm
    z_a \norm e_a\ \hbox{ with equality if, and only if,}\ 
  z=\gamma e,\ \hbox{for}\ \gamma \in \C.
  \end{equation}
  Since $e$ is maximal, $\norm e_a=\sqrt n$ and by \ref{algnorm} $\|z\|_a\leq \sqrt n\|z\|$ giving
     \begin{equation}\label{ineq}
     \norm
    z_a \norm e_a \leq n\ \hbox{for all}\  z \in \ol B.
      \end{equation}    
      Therefore 
        $\absval{\frac{\ip ze_a}{\ip ee_a}}
    \le 1$ for all $z\in \ol B$ with equality if, and only if, we have equality in both (\ref{CS}) and (\ref{ineq}), namely  
    \begin{equation}\label{seclast} 
       \absval{\frac{\ip ze_a}{\ip ee_a}}\leq 1 \ \hbox{with  equality if, and only if,}\ 
  z=\gamma e,\ \hbox{for}\ \gamma \in \mathbb{T}=\del \Delta. 
   \end{equation}
   Since for $ \mu \in    \mathbb{T},\  \frac{1}{2}|1+\mu|=1$
 precisely when $\ \mu=1$, 
      it follows  from \eqref{eq:1} and \eqref{seclast} that $\norm{h(z)}\leq 1\ \hbox{for all}\  z \in B, \ \hbox{with equality if, and only if,}\   z=e.$
\end{Proof}

In the absence of finite dimensionality, the above result still
allows that the function $h$ may have norm determined (though not attained) 
away from $e$.  To address this we require the
following Lemma.
  \begin{Lemma}\label{lem:elem}
    Let $H$ be an inner product space, $e, z \in H$ with $\norm z \le
    \norm e$ and $\epsilon > 0$.  There exists $\delta>0$ such that
    $\absval{1- \frac{\ip ze}{\ip ee}}<\delta$ implies $\norm{z-e}<\epsilon$.
  \end{Lemma}

  %{\color{blue} FOR REMOVAL- just here for proofreading.

    \begin{Proof}
      Write $z^\perp = z - \frac{\ip ze}{\ip ee}e$ so that $z=
      \frac{\ip ze}{\ip ee}e + z^\perp$ where $\ip{z^\perp}e =0$.
      Then
      \( \norm{z^\perp}^2+ \frac{\absval{{\ip ze}}^2}{\ip ee} = \norm z^2 \le \norm
        e^2\) which implies \[ \norm{z^\perp} \le \norm e\sqrt{ 1-
          \frac{\absval{\ip ze}^2}{\ip ee^2}}.\]
      Now $\absval{1- \frac{\ip ze}{\ip ee}}<\delta$ gives
      $1-\absval{\frac{\ip ze}{\ip ee}} <\delta$ and $1-\absval{\frac{\ip ze}{\ip ee}}^2 <2\delta$  so $\norm{z^\perp}
      < \norm e \sqrt{2\delta}$.  Finally
      \begin{align*}
        \norm{z-e} &= \norm { (z-z^\perp) -e + z^\perp} \\
                   &\le \norm{\frac{\ip ze}{\ip ee}e -e } + \norm{z^\perp}
        < \delta\norm e + \sqrt{2\delta} \norm e
      \end{align*}
      which is smaller than $\epsilon$ for $\delta>0$ sufficiently small.
      \end{Proof}
     % }
  
  \begin{Proposition}\label{prop:C}
    Let $Z$ be a finite rank \JBT\ and X be a Banach space.  For $e\in
    \Gamma$, there exists a holomorphic map $h:Z\to X$ such that 
    $\norm f_{\ol B} =1$ but $\norm f_{\ol B \setminus B(e,\epsilon)}
    <1$ for $\epsilon >0$.
  \end{Proposition}

    \begin{Proof}
    The same function $h$ provided in Proposition~\ref{prop:B}  suffices.  Observe 
    that for $\mu \in \ol\Delta=\ol \Delta(0,1)$ and $\eta\in(0,1)$ then
    $\half\absval{1+\mu}>1-\eta$ implies $\mu\in
    \ol\Delta(0,1)\setminus \Delta(-1,2-2\eta)$ and hence $\mu\in
    \Delta(1,\delta)$ where $\delta^2+ (2-2\eta)^2=2^2$, that is,
    $\delta=2\sqrt{2\eta-\eta^2}$.  

    Let $\epsilon>0$ and choose $\delta$ according to
    Lemma~\ref{lem:elem} where $\ip zz_a \le \ip ee_a =n$.
    Choose $\eta>0$ according to the observation above so that
    $\half\absval{1+\mu}>1-\eta$ implies $\mu \in \Delta(1,\delta)$.

    Now, suppose $\norm{h(z)}> 1- \eta$ so that $\half \absval{1+
      \frac{\ip ze_a}{\ip ee_a}}>1-\eta$.  Then $\absval{1-  \frac{\ip
        ze_a}{\ip ee_a}} <\delta$ and thus $\norm{z-e}_a<\epsilon$.
    As $\norm{z-e}\le \norm{z-e}_a$ we conclude
    $\norm{h}_{B\setminus B(e,\epsilon)} < 1-\eta<1$ as required. 
      
  \end{Proof}

We are now in a position to prove that the set of extreme points $\Gamma$ is the infinite dimensional analogue of the
Bergmann-Shilov boundary for finite rank JB*-triples.

\begin{Theorem}\label{bergshilov}
  Let $Z$ be a finite rank \JBT.  Then the set $\Gamma$ of extreme points of $\ol B$ is the
  smallest closed subset $\Lambda$ of $\ol B$ such
  that \begin{equation} \sup\{\norm{f(z)}:z\in\ol B\} = \sup \{
    \norm{f(z)}: z\in \Lambda\}\label{eq:det}\end{equation}
 for all $f:B \mapsto \C$ holomorphic on $B$ with continuous extension to $\del B$.
  \end{Theorem}

  \begin{Proof}
    Corollary~\ref{determining} shows that $\Gamma$ is a (closed) set
    satisfying~\eqref{eq:det}.
 Now let $\Lambda$ be any
    closed set in $\ol B$ satisfying~\eqref{eq:det}.  Suppose $e\in\Gamma$ but
    $e\notin \Lambda$.  Then as $\Lambda$ is closed, $B(e,\delta)\cap
    \Lambda=\emptyset$ for some $\delta>0$.   By
    Proposition~\ref{prop:C}, there exists $h\in \H(Z, \C)$ such
    that $\norm h_{\Lambda}  \le \norm h_{\ol B \setminus B(e,\delta)} < \norm h_{\ol B}$, which is a
    contradiction to \eqref{eq:det}.  We conclude $e\in \Lambda$ and
    consequently, $\Gamma \subset \Lambda$ as required.
  \end{Proof}

  If a finite rank \JBT\  $Z$ contains a unitary tripotent then
  $\Gamma_1$ is a non-empty closed subset of $\Gamma$ which, by
  Proposition~\ref{determining} part (2), satisfies \eqref{eq:det}.  By
  Theorem~\ref{bergshilov} then, $\Gamma_1$ cannot be a proper subset
  of $\Gamma$ and we have the following consequence.
  
  \begin{Corollary}\label{maximalisunitary}
    If a finite rank \JBT\ has a unitary tripotent then all of its
    maximal tripotents are unitary.
  \end{Corollary}
  
  \begin{Corollary}\label{last}
    Let $Z$ be a finite rank \JBT\ with open unit ball $B$ and $v\in
    \del B$. Then 
    $\Gamma\subseteq \ol {G_v^0}$.    
  \end{Corollary}

  \begin{Proof}
    From Proposition~\ref{vorbits}, $\ol{G_v^0}$ is a closed
    determining set satisfying \eqref{eq:det} so the result follows by Theorem~\ref{bergshilov}.
  \end{Proof}

\bibliographystyle{acm}
%\bibliography{bib.bib}
\def\cprime{$'$}

\noindent
%{P. Mellon,}\newline
%{School of Mathematics,}\newline
%{University College Dublin}, \newline
%{Dublin 4, Ireland.} \newline

%email: {mackey@maths.ucd.ie, pmellon@maths.ucd.ie}

\end{document}